\documentclass[12pt, reqno]{amsart}

\usepackage{amsthm, amsmath, amssymb}
\usepackage{color,epsfig, graphicx}
\usepackage{url}

\begin{document}

\title{Socially Optimal Charging Strategies for Electric Vehicles}
\author{Elena Yudovina~and~George Michailidis}
\thanks{Department of Statistics, University of Michigan, Ann Arbor, MI, 48109 USA}
\thanks{The research of the first author was supported by NSF grant DMS-1204311. The research of the second author was supporter by NSF grant DMS-1106695 and NSA grant H98230-10-1-0203.}
\maketitle

\newtheorem{theorem}{Theorem}
\newtheorem{lemma}[theorem]{Lemma}
\newtheorem{assumption}{Assumption}


\newcommand{\BR}{{\mathbb R}}
\newcommand{\BN}{{\mathbb N}}

\newcommand{\BE}{{\mathbb E}}
\newcommand{\BP}{{\mathbb P}}

\newcommand{\CM}{{\mathcal M}}
\newcommand{\CL}{{\mathcal L}}
\newcommand{\CS}{{\mathcal S}}
\newcommand{\CC}{{\mathcal C}}
\newcommand{\DD}{{\mathcal D}}
\newcommand{\CJ}{{\mathcal J}}
\newcommand{\CI}{{\mathcal I}}
\newcommand{\CE}{{\mathcal E}}

\newcommand{\one}{{\mathbf 1}}

\newcommand{\ov}{\overline}
\newcommand{\ul}{\underline}
\newcommand{\del}{\partial}
\newcommand{\trace}{\text{trace}\,}

\newcommand{\weakto}{\stackrel{w}{\rightarrow}}

\providecommand{\abs}[1]{\left\lvert#1\right\rvert}
\providecommand{\norm}[1]{\lVert#1\rVert}

\newcommand{\floor}[1]{\lfloor#1\rfloor}
\newcommand{\ceiling}[1]{\lceil#1\rceil}

\newcommand{\beql}[1]{\begin{equation}\label{#1}}
\newcommand{\eqn}[1]{(\ref{#1})}

\renewcommand{\labelenumi}{\theenumi.}

\begin{abstract}
Electric vehicles represent a promising technology for reducing emissions and dependence on fossil fuels and have started entering different automotive markets. In order to bolster their adoption by consumers and hence enhance their penetration rate, a charging station infrastructure needs to be deployed. This paper studies {\em decentralized} policies that assign electric vehicles to a network of charging stations with the goal to achieve little to no queueing. This objective is especially important for electric vehicles, whose
charging times are fairly long. The social optimality of the proposed policies is established in the many-server regime, where each station is equipped with multiple charging slots. Further, convergence issues of the algorithm that achieves the optimal policy are examined. Finally, the results provide insight on how to address questions related to the optimal location deployment of the infrastructure.
\end{abstract}


\section{Introduction and model}
There has been an increasing penetration of Plug-in Hybrid and pure Electric Vehicles (PHEV/EV) over the least few years
\cite{IEA} . This is due to developments in battery technology that have dramatically increased their range \cite{Lukic_et_al}, advances in charging technologies that have reduced their charging times \cite{Price}, incentives that have lowered their acquisition and operation cost and an overall desire to lower emissions \cite{MIT_report}.

On the other hand, they represent a potential source of disruption to normal grid operations if not integrated carefully because they need to connect to the distribution network to charge \cite{Boulanger_et_al, Hajimiragha_et_al}. Currently, EVs are equipped primarily with lithium-ion batteries, ranging in energy capacity from 5 kWh for short-range PHEVs to 50 kWh for high performance EVs. Further, today's and future EVs are designed with a wide range of specifications to satisfy different customer preferences.

Their impact on power grid operations
will heavily depend on their market penetration \cite{Price}. Estimates vary widely, ranging from 3 to 18 million vehicles by 2025 and from 5 all the way to 40 million vehicles (approximately 20\% of the total US market) by 2030 \cite{MIT_report}. Their disruptive impact is mainly due to the energy load they represent. On average, under normal charging conditions (1.4 kW) an EV represents a 1.3-fold of a full household load, whereas fast charging technologies (7.2 kW) correspond to an almost 3-fold increase \cite{MIT_report, Price}. The injection of such large loads coupled with the possible uneven geographic distribution of EVs would definitely strain the entire grid as argued in \cite{Galus_et_al}. 

Given the current predominance of PHEV vehicles, the literature has largely focused on scheduling at home overnight charging (see \cite{Ma_et_al, Pang_et_al, Verjilbergh_et_al} and references therein). The proposed approaches treat the induced load as an aggregate and discuss different mechanisms on how to shift it during night hours to take advantage of the underutilized electricity generation assets. However, with increasing penetration rates, efficient operation of an expanding charging station infrastructure becomes a key issue. Charging EVs is a rather slow process, as even fast charging takes at least half an hour \cite{MIT_report}, thus requiring careful scheduling policies to provide the necessary quality of service to customers. Faster charging technologies (e.g. DC charging) could mitigate some of these effects, but as mentioned above, electric utilities have concerns about possible negative impacts of such technologies on the power grid, if deployed at large scale. The work to date on charging stations has mostly focused on modeling and optimizing the architecture of a {\em single} charging station \cite{Bai_et_al, Bayram_et_al, Reed_et_al}. 

In this paper, we consider {\em decentralized dispatching policies} that assign electric vehicles to a {\em network} of charging stations.
Our focus is on guaranteeing quality-of-service to charging customers. Due to the lengthy charging times, our primary focus is to ensure little to no customer queueing. Our results apply to the ``many-server'' regime, when each charging station is equipped with multiple charging slots (e.g. a commercial parking lot, a strip of street parking in a
densely urban area), possibly in a mix of charging technologies (e.g. Level-2 charging, coupled with very fast DC charging infrastructure).

We provide next a description of the modeling framework.
Consider  EVs and a network of charging stations within a specific geographical region, 
for example an urban area or a section of the highway system. A vehicle that needs to recharge broadcasts a signal indicating its location and battery type and status, and receives responses from the charging stations in its neighborhood. We assume that the EV has preferences among the different charging stations, which we model through costs; these 
may be related to the distance of the current location of the EV to the charging station,  or simply whether it is capable of reaching the station using its remaining battery power. 
Based on the stations' responses, the EV chooses a charging station and immediately proceeds there. We are interested in designing {\em socially optimal} charging strategies; namely, to  engineer the response signals from stations to vehicles and the resulting {\em routing decisions}, so that the average cost incurred by the vehicles is minimized.

We will make the simplifying assumption that the vehicles can be partitioned into a finite set of types $i=1,\dotsc,I$ which encode their location, preferences, and battery technology. Thus, the signal that the EV broadcasts to the nearby charging stations indicates the type of the vehicle. We assume also that the various charging slots available in the
network can be partitioned into a finite set of types $j=1,\dotsc,J$, encoding their location and technology. We are interested in the regime where there are many charging slots of each type: say $N_j$ chargers of type $j$, with $N_j$ large (e.g. at least 10). Thus, as mentioned above, the architecture of the charging station in the network is such that it comprises
of many identical chargers.

Let $\lambda_i$ be the rate at which EVs of type $i$ make charging requests. Let $\mu^{-1}_{ij}$ be the expected time it takes for a charger of type $j$ to satisfy a request of type $i$; this may be infinite (corresponding to $\mu_{ij}=0$), e.g. if an EV cannot reach the charger or their technologies are incompatible. Finally, let $c_i(j)$ be the cost that an EV of type $i$ incurs on being assigned to type $j$; this may be a measure of the distance between the vehicle and the chargers of type $j$, or a constant whenever the vehicle can reach the charger. We set $c_i(j) = \infty$ whenever an EV of type $i$ can not reach a charger of type $j$. Without any guidance from the system/network, we would expect a vehicle of 
type $i$ to always choose a charger of type $j$ that minimizes $c_i(j)$; however, the following example shows this behavior to be suboptimal in some cases.

Consider a system with two vehicle classes $A$ and $B$, and two charger classes $1$ and $2$. Suppose that $A$ cannot reach class $2$, and $B$ can use both charger classes, but prefers type $1$. That is, $c_A(2) = \infty$ and $c_B(1) < c_B(2)$. In this case, if both EV types preferentially go to chargers of type $1$, vehicles of type $A$ may have to queue. If, on the other hand, vehicles of type $B$ can be directed to chargers of type $2$ when few chargers of type $1$ remain, queueing can be avoided.

Another example is obtained by considering two EV and charger classes as before. Suppose that $A$ prefers type $1$ and $B$ prefers type $2$; that is, $c_A(1) < c_A(2)$ and $c_B(2) < c_B(1)$. Suppose further that $\mu_{A1} < \mu_{A2}$ while $\mu_{B1} = \mu_{B2}$, meaning that it is faster to charge vehicles of type $A$ at the more distant charger. Clearly, in heavy traffic we must encourage vehicles of type $A$ to charge at the more distant, fast charging station $2$, while vehicles of type $B$ will correspondingly need to be directed to the more distant (and not even faster!) charging station $1$.

From a social perspective, the main objective is to to design charging strategies for a network of stations, where such routing decisions will happen automatically.

This problem has similarities to the inventory and facility location problem: we are interested in ``distributing'' a finite supply of chargers among vehicles, subject to location constraints, in such a way as to keep vehicles from being ``undersupplied''. Extensive literature exists on such problems; see, for example, \cite{ORTextbook} and references therein. However, there are two key differences in our set-up. The first is that EV demand is mobile. Thus, demand from the vehicles of type (location) $i$ can be split between several different charging stations (facilities). In fact, our work shows that this splitting should be encouraged. The implication from an algorithmic perspective is  that instead of solving an integer program that typically the case in the standard facility location problem, we need to solve its convex relaxation. The second key difference between our setting and the inventory / facility location problem is that we wish to avoid centralized decisions. Instead, for our set-up to be scalable to the size of a large urban or even a metropolitan area with hundreds of charging stations and thousands of vehicles, we need the decision-making to be distributed. The latter goal justifies our formulation for letting the EVs 
 choose their preferred charging station based on the information communicated to them by nearby stations, rather than involving some centralized planning scheme.

To achieve this goal and design the required efficient, {\em distributed} algorithm for routing EVs to charging stations so as to avoid excessive delays due to queueing, 
we employ ideas from queueing theory and communication systems. We introduce the GPD algorithm (Greedy Primal-Dual) that has been successfully used in the call
center literature and establish analytic guarantees for its performance in large scale network of charging stations. Its main feature is its distributed an online nature,
and also its automatic adaptability to changing arrival patterns. We also present two variants (Load Balancing and Freest Charger Shortest Queue), which are shown to 
exhibit superior performance to the GPD algorithm in selected settings.

\subsection{Modeling Assumptions}
We discuss next the main modeling assumptions. First, we have stochastic assumptions on key processes of the system under consideration,
that lead to the analytic guarantees on the system behavior. Second, we specify which of the system parameters are known to which participants; the latter 
assumptions ensure that the system behaves in a distributed manner.

We begin with stochastic process assumptions.
\begin{assumption}\label{ass: FLLN}
The arrival process of charging requests of EVs of type $i$ satisfies a functional law of large numbers approximation.
Let $A_i(t) = \#\{\text{requests of type $i$ up to time $t$}\}$; we require
\[
\frac1r A_i(rt) \implies \lambda_i t, ~~ \text{as $r \to \infty$,} ~~ \text{uniformly on compact sets.}
\]
Throughout, the notation $\implies$ incicates uniform convergence on compact sets. Furthermore, we require that arrivals have bounded second moment: $\BE[\bigl(A_i(t+1) - A_i(t)\bigr)^2] < \infty$ uniformly in $t$.

The process of service completions of EVs of type $i$ by chargers of type $j$ satisfies a functional law of large numbers approximation. Specifically, let $S_{ij}(t)$ be the number of EVs of type $i$ that have completed service with one charger of type $j$ when that charger has spent a total amount of time $t$ charging vehicles of type $i$; we require
\[
\frac1r S_{ij}(rt) \implies \mu_{ij} t, ~~ \text{as $r \to \infty$,} ~~ \text{u.o.c.}
\]
We also require that service times have bounded second moment: $\BE[\bigl(S_{ij}(t+1) - S_{ij}(t)\bigr)^2] < \infty$ uniformly in $t$.
\end{assumption}
These are quite general assumptions, that are satisfied for example when interarrival and service times are independent and identically distributed (iid) and possessing
 finite second moments.

Assumption~\ref{ass: FLLN} will be taken to hold throughout the paper. It is sufficient to show that the proposed algorithm has optimal throughput, meaning that it will successfully recharge all EVs whenever it is possible to do so (perhaps with queueing delays). However, in order to provide more precise guarantees, for example on the probability of an arriving vehicle finding a free charger, we will need more control over the deviations of the arrival and service processes from their FLLN approximation. 
Assumption~\ref{ass: FCLT} asserts that the arrival and service processes obey a functional central limit theorem.
\begin{assumption}\label{ass: FCLT}
The arrival processes of requests of type $i$ obey the functional central limit theorem:
\[
\frac{1}{\sqrt{r}}\bigl(A_i(rt) - \lambda_i r t\bigr) \implies W_i(t), ~~ \text{as $r \to \infty$},~~ \text{u.o.c.}
\]
Here, $W_i$ are a set of independent Brownian motions (one per EV type) with some finite variance.

In addition, the service processes of EVs of type $i$ by chargers of type $j$ obey the functional central limit theorem:
\[
\frac{1}{\sqrt{r}}\bigl(S_{ij}(rt) - \mu_{ij} r t\bigr) \implies W_{ij}(t), ~~ \text{as $r \to \infty$},~~ \text{u.o.c.}
\]
Here, $W_{ij}$ are a set of independent Brownian motions (one per pair of vehicle type and charger type) with some finite variance.
\end{assumption}

These assumptions are standard in the queueing literature. They hold when the interarrival and service times are iid with finite second moment, and typically allow 
one to model the queueing process as a reflected Brownian motion. We will comment further on this in Section~\ref{Halfin-Whitt}.

We next specify what information is available to which participants. The system has parameters $\lambda_i$ (arrival rate of requests of type $i$), $\mu_{ij}$ (service rate of requests of type $i$ by chargers of type $j$), $N_j$ (number of chargers of type $j$), and $c_i(j)$ (the cost associated with assigning a vehicle of type $i$ to a charger of type $j$). In addition, the scheduling algorithm will use a parameter $\beta$.
\begin{assumption}\label{ass: who knows what}
The number of possible types is finite.\\
The parameters $\lambda_i$ are unknown to anyone in the system; the algorithm will implicitly estimate them.\\
The parameters $\mu_{ij}$ are assumed to be known to a request for which they are relevant (e.g., included as part of the exchange between the EV 
and the charging stations); if $\mu_{ij}$ is not communicated, it is assumed to be 0.\\
The costs $c_i(j)$ are assumed to be known to the EV for which it is relevant. Specifically, we require the EV to be able to compare quantities of the form $c_i(j) + K_{ij} \mu_{ij}^{-1}$, where $K_{ij}$ will be a quantity communicated by the charging station.\\
The parameter $\beta$ (a small real number) is assumed to be the same at all charging stations.
\end{assumption}
Note that all communications and knowledge are {\em local}: EVs only acquire information about nearby stations, and stations only acquire information about the nearby requests.

The remainder of the paper is organized as follows:
In Section~\ref{sec: policies} , we formulate the main algorithm, GPD, and show its throughput optimality. In Section~\ref{sec: limiting regimes}, 
we present a detailed analysis of the behavior of large systems, provide Brownian motion approximations for queue sizes and the associated probability of finding a free charger upon arrival. We also state a key result (Theorem~\ref{thm: FCLT},) whose proof is given in the Appendix. In Section~\ref{sec: LB}, we introduce the Load Balancing (LB) algorithm, which is designed to spread the load more evenly between chargers to avoid excess queueing. Section~\ref{sec: 0-infty} discusses the case when the costs $c_i(j) \in \{0,\infty\}$, i.e. EVs are indifferent between chargers provided they are within reach. For the case of no user costs, we present another algorithm, Freest Charger Shortest Queue (FCSQ), which reacts faster  than GPD or LB to changes in the arrival pattern. Section~\ref{sec: simulations} presents selected simulation results of the behavior of the three algorithms on a simple system. 
Finally, Section~\ref{sec: discussion} discussed insights into planning the charging network, as well as directions for future research. 

\section{Scheduling policies and throughput optimality}\label{sec: policies}

In this section, we formulate the policy that will achieve socially optimal average costs while keeping queueing low, if possible. We begin by formulating the corresponding linear program. Let $\lambda_{ij}$ be the average rate at which vehicles of type $i$ are routed to stations of type $j$. The basic linear program we consider is as follows.
\begin{subequations}\label{STABILITY}
\begin{alignat}{2}
&\text{minimize } && \sum_{i,j} \lambda_{ij} c_i(j)\\
&\text{s.t. } && \lambda_i = \sum_j \lambda_{ij},~~~~ \forall i\label{routing}\\
&&& N_j \geq \sum_i \frac{\lambda_{ij}}{\mu_{ij}}, ~~~~ \forall j\label{capacity}\\
&\text{over } && \lambda_{ij} \geq 0, ~~~~ \forall i,\,j.
\end{alignat}
\end{subequations}
The objective function here is the rate at which costs are incurred. The constraints are the basic feasibility ones: all arriving requests need to be assigned, but on average no more than $N_j$ chargers of type $j$ may be used. Note that the optimal solution may be infinite, corresponding to insufficient capacity in the system. Call the value of this linear program $S(\lambda)$.

We define the feasible region to be the set of routing rates for which the solution to this linear program is finite and the capacity constraints \eqref{capacity} are strictly satisfied. 
That is,
\[
\Lambda = \{(\lambda_i)_{i=1,\dotsc,I}:~ S(\lambda) < \infty, \text{ and } N_j > \sum_j \frac{\lambda^*_{ij}}{\mu_{ij}}~ \forall j\}
\]
for some optimal solution $\lambda^*_{ij}$ of \eqref{STABILITY}.

We now introduce the Greedy Primal-Dual (GPD) Algorithm of \cite{GPD}, which will determine a scheduling rule that implicitly solves the linear program \eqref{STABILITY}. The parameter $\beta$ is a small real number that determines the trade-off between convergence speed and precision of the solution. While the algorithm is a special case of the one in \cite{GPD}, technical differences between our setting and that in \cite{GPD} imply that analytic guarantees presented in Section~\ref{Halfin-Whitt} are not automatic and
need to be established rigorously.

\vspace{1cm}\noindent{\bf GPD Algorithm:}
\begin{enumerate}
\item Each charger type maintains a virtual queue variable $Q_j(t)$ (this need not be an integer). The latter are appropriately initialized, e.g. $Q_j(0) = 0$.
\item When an EV of type $i$ requests service at time $t$, we locate
\[
j^* \in \arg\min_j c_i(j) + \frac{\beta Q_j(t)}{\mu_{ij}}.
\]
That is, all stations in the neighborhood of $i$ communicate parameters $\beta Q_j(t)$ (or $\beta Q_j(t) / \mu_{ij}$), and the user picks the station $j^*$ that minimizes $c_i(j) + \beta Q_j(t) / \mu_{ij}$. The vehicle announces its decision to $j^*$, and the corresponding virtual queue is incremented:
\[
Q_j^* \mapsto Q_j^* + \frac{1}{\mu_{ij^*}}.
\]
\item Decrease \emph{all} virtual queues at rate $N_j$ per time unit, provided they are positive. (Once a virtual queue hits 0, it stays at level 0 until some EV is routed to it.)
\end{enumerate}
Note that the algorithm runs in continuous time and {\em no synchronization} is necessary between different stations.

One important instance of the algorithm is when $c_i(j) \in \{0,\infty\}$ for all $i,j$; that is, the EVs are indifferent among the charging stations as long as 
they can use them at all. We will discuss this case in detail in Section~\ref{sec: 0-infty}. For now we point out that in this case, the parameter $\beta$ is unnecessary.

The results of \cite{GPD} imply that whenever $\Lambda \neq \emptyset$, the GPD algorithm will stabilize all virtual queues. This intuitively means that the algorithm will be routing vehicles to chargers at rates $\lambda_{ij} \in \Lambda$. The more precise statement is as follows. Let $A_{ij}(t)$ denote the number of vehicles of type $i$ that have been routed to chargers of type $j$ up to time $t$. The limit $\lim_{t \to \infty} \frac1t A_{ij}(t)$, denoting the average rate, may not exist; however, any sequence of times $t \to \infty$ has a subsequence along which such the limit does exist. Pick one such set of subsequential limits, and call it $(\lambda_{ij})$. Then, if $\Lambda \neq \emptyset$, the rates $\lambda_{ij}$ satisfy the routing constraint \eqref{routing} and the capacity constraint \eqref{capacity}. In general, we may have $\sum_{ij} \lambda_{ij} c_i(j) > S(\lambda)$, i.e. the rates may not be socially optimal; note, however, that the algorithm remains throughput-optimal, that is, it stabilizes the system whenever the arrival rates belong to $\Lambda$. In Section~\ref{sec: limiting regimes} we will see that as $\beta \to 0$, the algorithm converges to the socially optimal routing rates.

Note that in the case of Poisson arrival and service process, and with service rates $\mu_{ij}$ being rational numbers,  the virtual queueing system is a 
countable state-space Markov process. In that case, there will be well-defined steady-state rates $\lambda_{ij}$ to which the above observations apply.

Further, note the natural form (from the user's point of view) of the routing decision: the EV driver is asked to add to her intrinsic costs $c_i(j)$ a certain charge per unit time $\beta Q_j(t)$. This charge will be greater for the stations that are in high demand, and lower for the stations that are less congested. We point out that, because the GPD algorithm uses virtual queues rather than actual queues to make routing decisions, it can be run in the background of some other scheduling mechanism, provided a model for the costs $c_i(\cdot)$ is available.

\section{Limiting regimes}\label{sec: limiting regimes}
We now describe the asymptotic behavior of the system in certain limiting regimes. Our interest is in large systems; consequently, we will be interested in the case $\lambda_i \to \infty$. To accommodate the increasing arrival rate, we will consider $N_j \to \infty$ (many chargers at each station), holding the number of charger types and the speed of charging fixed. We will also consider the effect of taking $\beta \to 0$, where $\beta$ is the parameter used in defining the GPD algorithm.

In order to state the results, we make an assumption on the solution structure of the linear program \eqref{STABILITY}. Some of the conclusions apply even when these do not hold; however, the exposition would become more cumbersome.

\begin{assumption}\label{ass: CRP}
The optimal solution $(\lambda^*_{ij})$ to the linear program is unique. The optimal set of dual variables $(q^*_j)$ is unique. The pairs $(ij)$ for which $\lambda^*_{ij} > 0$ are called \emph{basic activities}.
\end{assumption}
This is a part of the \emph{complete resource pooling} assumption, which is commonly made in many-server queueing literature, e.g. \cite{Shadow_routing}. As pointed out in \cite{Shadow_routing}, Assumption~\ref{ass: CRP} holds for generic values of parameters $\lambda_{ij}$, $\mu_{ij}$, and $c_i(j)$.

We describe briefly the techniques used to obtain the stability and ``fluid-scaled'' convergence results in Sections~\ref{sec: beta fixed}--\ref{sec: N_j fixed}, relying only on Assumption~\ref{ass: FLLN} for the underlying stochastic processes. The technique is standard in the queueing literature, and involves the use of fluid limits rather than constructions of explicit Lyapunov functions. Specifically, one proves a functional law of large numbers approximation (``fluid limit'') for the stochastic processes involved. If the resulting limiting trajectories satisfy certain properties, then one can conclude that the original stochastic system is stable. (For a system described by a Markov process, by ``stable'' we mean positive Harris recurrent.) A good exposition can be found in \cite{Bramson}. This technique is often extended by using a Brownian approximation for the underlying stochastic processes of arrivals and service completions. Such assumptions allow a description of the queue lengths on a finer (``diffusion'') scale, where they are non-vanishing; specifically, queues are approximated by a reflected Brownian motion. These techniques present greater technical challenges, particularly in deriving steady-state results. We discuss diffusion-scaled approximations in Section~\ref{Halfin-Whitt}, where our exposition follows closely the similar results of \cite{Shadow_routing}. A good reference on diffusion approximations of queueing networks is \cite{Harrison_Williams}.

Next, we examine specific cases.

\subsection{$N_j \to \infty$, $\beta$ fixed}\label{sec: beta fixed}
We begin by considering the regime in which $\beta$ is fixed, but the arrival rates and the number of chargers tend to infinity. More formally, we consider a sequence of systems indexed by $r$, with $\lambda^r_i = r \lambda_i$ and $N^r_j = r N_j$. (Note that $^r$ is a superscript, not the $r^{\text{th}}$ power.) In this regime, by results from \cite{GPD} which we already stated above, the average routing rates $\lambda^r_{ij}$ are guaranteed to stabilize the system provided $(\lambda_i) \in \Lambda$, but are not in general guaranteed to be socially optimal even as $r \to \infty$. One important result comes from considering large systems. If $(\lambda_i) \in \Lambda$, the system is strictly underloaded. Consequently, as $r \to \infty$, the probability of an arrival request having to queue and the average queueing time will both converge to zero. This is an important consideration in the problem of charging electric vehicles.

\subsection{$N_j$ fixed, $\beta \to 0$}\label{sec: N_j fixed}
This is the regime in which the GPD algorithm converges to the optimal solution of the linear program \eqref{STABILITY}. We consider a sequence of systems indexed by $r$, for which $\beta^r \to 0$ and all other parameters are unchanged. This corresponds to running a sequence of slightly different GPD algorithms on the same external pattern of requests.

All the results for this case, follow straightforwardly from \cite{GPD}; we rewrite them in the notation of our problem. Let $Q^r_j(t)$ be the sequence of virtual queues. Then, uniformly over any compact time interval of $t \in [0,T]$, we have
\begin{subequations}\label{eqn:GPD}
\begin{equation}\label{Q uoc}
\beta^r Q^r_j(\beta^r \cdot) \implies q_j(\cdot), ~~~\text{u.o.c. as} ~ \beta^r \to 0,
\end{equation}
where the limiting trajectories satisfy
\begin{equation}\label{q to q*}
q_j(t) \to q_j^*, ~~ \text{as} ~ t \to \infty.
\end{equation}
\end{subequations}
Here, the parameters $(q_j^*)$ are the optimal dual variables corresponding to the capacity constraints \eqref{capacity}.

In addition, we describe the convergence of the routing rates. Let $X^r_{ij}(t)$ be the exponentially-weighted average of the routing decisions up to time $t$. Then, uniformly on compact sets,
\begin{subequations}
\begin{equation}\label{X uoc}
X^r_{ij}(\beta^r \cdot) \implies x_{ij}(\cdot), ~~~\text{u.o.c. as} ~ \beta^r \to 0;
\end{equation}
where the limiting trajectories satisfy
\begin{equation}\label{x to x*}
x_{ij}(t) \to \lambda^*_{ij}, ~~ \text{as} ~ t \to \infty.
\end{equation}
\end{subequations}
Here, $\lambda^*_{ij}$ are the (unique) optimal solution of the linear program \eqref{STABILITY}.

These results only require the functional law of large numbers given in Assumption~\ref{ass: FLLN}, but not the functional central limit theorem, as shown in  
 \cite{GPD}. If the interarrival and service times are assumed to be Poisson, and the parameters $\mu_{ij}$ are rational (so that the virtual queueing system is a countable state-space Markov process), these results imply that the steady-state routing rates $\lambda^r_{ij}$ converge to the optimal rates $\lambda^*_{ij}$ as $\beta^r \to 0$.

\subsection{$N_j \to \infty$, $\beta \to 0$}\label{Halfin-Whitt}
In this case, we combine the effects of the previous two limits. The operational regime we are interested in is known as the Halfin-Whitt regime \cite{Halfin_Whitt}. When charging stations have many individual chargers, it is possible to operate a heavily loaded system while keeping waiting times short and providing a service guarantee on the probability of an arriving request having to wait. This effect was originally observed in a queueing system with a single pool of many servers, but has since been shown to apply in settings with multiple server types, at least under conditions similar to Assumption~\ref{ass: CRP}. This result leads to the ``square-root staffing'' principle observed in \cite{Halfin_Whitt}.

From the previous sections, we expect to see virtual queues $\beta Q_j \approx q^*_j$. Our interest will be in the deviations from these values. If $\beta$ is held fixed, these deviations may, in general, be large. We will now choose $\beta \to 0$ along an appropriate sequence as the size of the system grows larger. We will establish that after rescaling time, the process $\beta Q(\cdot)$ can be described as a Brownian oscillation around $q^*$. Up to a multiplicative constant, this result is the best possible, whenever the underlying arrival and service completion processes satisfy the functional central limit theorem given in Assumption~\ref{ass: FCLT}. It can be translated into a reflected Brownian motion approximation for the actual queue sizes in the system.

We consider a sequence of systems indexed by $r \to \infty$. The request arrival rates satisfy $\lambda^r_i = r \lambda_i$. The sizes of the charger pools satisfy $N^r_j = r N_j + 
\sqrt{r} n_j + O(1)$, for some collection of values of $n_j \in \BR$ (possibly negative); this is the ``square-root staffing''. (The $O(1)$ term is included to ensure that $N^r_j$ is an integer.) We assume that the physical limits on the charging rate to not change, so that the charging rates $\mu_{ij}$ do not change with $r$. Finally, we choose $\beta^r = f(r)^{-1}$ for some function $f(r)$ with $r^{1/2} \ll f(r) \ll r$; for example, we may take $f(r) = r^{3/4}$. (Larger values of $f(r)$, corresponding to smaller values of $\beta^r$, will lead to more precise but slower convergence.)

Let $A^r_{ij}(t)$ be the number of EVs of type $i$ routed to chargers of type $j$ during the interval $[0,t]$. In Theorem~\ref{thm: FCLT}, we will establish that this quantity can be approximated by $\lambda^r_{ij} t + \sqrt{r} B(t)$ for some Brownian motion $B$, provided we consider the system at a time when the GPD algorithm has reached its steady-state. The 
result is similar in nature to  \cite{Shadow_routing}, but there are technical differences discussed in the proof.

We state the next result next, whose proof is given in the Appendix.
\begin{theorem}\label{thm: FCLT}
Assume the arrival process satisfies the FCLT given in Assumption~\ref{ass: FCLT}. Let
\[
\hat q^r_j(t) = r^{-1/2} \left(Q^r_j(t) - (\beta^r)^{-1} q^*_j \right),
\]
and let
\[
\hat a^r_{ij}(t) = r^{-1/2} \left(A^r_{ij}(t) - \lambda^*_{ij} t \right).
\]
If $(\hat q^r_j(0)) \to 0 \in \BR^J$, then $\hat a^r_{ij}(\cdot) \implies H(W)$, where $W$ is the Brownian motion identified in Assumption~\ref{ass: FCLT}, and $H$ is a linear mapping defined in \eqref{eqn: H^{ST}} below. (Thus, $H(W)$ is also a Brownian motion, but with correlated components.)

Further, suppose the arrival and service completion processes are Poisson, and the parameters $\mu_{ij}$ are rational, so that the virtual queueing system becomes a countable state-space Markov process. For each $r$ consider the associated stationary version of the process. Then,  $(\hat q^r_j) \to 0 \in \BR^J$.
\end{theorem}

In order to define the load-balancing linear map $H$, we introduce some additional terminology. We call the activities $(ij)$ with $\lambda_{ij} > 0$ \emph{basic} and let $\CE$ be the set of basic activities. By Assumption~\ref{ass: CRP}, the solution to \eqref{STABILITY} is unique. It follows that the bipartite graph with vertices \{car types, charger types\} and (undirected) edges corresponding to basic activities is acyclic; that is, it is either a tree or a union of trees. We can now define the linear map $H$ as follows. For a vector $v = (v_1, \dotsc, v_I) \in \BR^I$, the image $w = H(v)$ (with coordinates indexed by basic activities $(ij)$) satisfies
\begin{align}\label{eqn: H^{ST}}
&\sum_j w_{ij} = v_i, ~~ \forall i\notag\\
&\frac{\sum_{i'} w_{i'j}/\mu_{i'j}}{\mu_{ij}} = \frac{\sum_{i''} w_{i''j'}}{\mu_{ij'}}, ~~ \forall i, (ij), (ij').
\end{align}
To see that this is a well-defined map, we show how to solve the above system of equations. Pick a leaf of one of the connected components of the basic activity graph, i.e. a vertex with a single edge coming out of it. On this edge, we either have $w_{ij} = v_i$, or can eliminate the variable $w_{ij}$ using the second set of equations. Repeating this process, we will arrive at a unique solution.

Next, we discuss the implications of the main result. The conclusion of the theorem asserts that the arrival process to each charger type can be described as a Brownian oscillation around the optimal rate $\sum_i \lambda^*_{ij}$; thus, the arrival process to each charger type satisfies a functional central limit theorem. Because we have assumed (Assumption~\ref{ass: FCLT}) that the service completion process satisfies it as well, we may approximate the number of occupied chargers by a Brownian oscillation about its optimal point, and the queue size as a reflected Brownian motion with drift; the drift is given by the collection of parameters $n_j$. In particular, the probability of an arriving request being asked to queue will depend on the quantity $\sum q^*_j n_j$, the sum being taken over the connected component of the basic activity tree containing the corresponding request type. (It will also depend on the linear map $H$ of \eqref{eqn: H^{ST}}.) Note that we require $O(\sqrt{r})$ overstaffing somewhere in the connected component of each EV type, but do not prescribe where. Consequently, it is advantageous to arrange the system to have large connected components to allow greater freedom in the placement of extra chargers.

\section{Load balancing}\label{sec: LB}

The GPD algorithm will typically produce rates which place a high load on one or more of the ``least expensive'' charger pools. This is reasonable in a large system, where square-root overstaffing is sufficient to deal with the high load, but may be undesirable in a system where the charging stations are small. One possibility for mitigating this is to reduce the value of $N_j$ in the GPD algorithm, so that on average only a fraction of the chargers may be used. However, this reduces the stability region of the algorithm. A better way to mitigate this difficulty is to encourage the charging stations to spread the load more evenly. Group the charging stations into \emph{clusters}. (Clusters may overlap.) Change the objective function of \eqref{STABILITY} to
\begin{subequations}\label{eqn:LB}
\begin{alignat}{2}
&\text{minimize } && \sum_{i,j} \lambda_{ij} c_i(j) + \sum_{l} W_l \rho_l\\
&\text{s.t. } && \lambda_i = \sum_j \lambda_{ij},~~~~ \forall i\label{routing2}\\
&&& N_j \geq \sum_i \frac{\lambda_{ij}}{\mu_{ij}}, ~~~~ \forall j\label{capacity2}\\
&&& \rho_l \geq \sum_i \frac{\lambda_{ij}}{N_j \mu_{ij}}, ~~~~ \forall j \in l\label{load}\\
&\text{over } && \lambda_{ij} \geq 0, ~~~~ \forall i,\,j.
\end{alignat}
\end{subequations}
where $l$ runs over the clusters, $W_l$ are weights, and $\rho_l$ is the maximal load of any charging station in the cluster. Modifying the objective in this manner means that charging stations within each cluster will ``try'' to have equal loads, because only the maximal load within a cluster is penalized. By adjusting the weight vector $W_l$, we can change the relative importance of spreading the load across different stations, and finding the lowest-cost routing pattern.

The corresponding modification of the GPD algorithm is called the Load Balancing (LB) algorithm.

\vspace{1cm}\noindent{\bf LB Algorithm:}
\begin{enumerate}
\item Each charger type maintains a virtual queue variable $Q_j(t)$, initialized e.g. to $Q_j(t) = 0$. Each charger type also maintains a virtual queue $L_j(t)$ for each cluster $L$ to which $j$ belongs. These are initialized so that $\sum_{j \in l} L_j(t) = W_l$ for each $l$. (We assume that this is possible.)
\item When a vehicle of type $i$ requests service at time $t$, we locate
\[
j^* \in \arg\min_j c_i(j) + \frac{\beta (Q_j(t) + L_j(t))}{\mu_{ij}}.
\]
The vehicle announces its decision to $j^*$, and the corresponding virtual queues are incremented:
\begin{align*}
Q_j^* \mapsto Q_j^* + \frac{1}{\mu_{ij^*}}, ~~ L_j^* \mapsto L_j^* + \frac{1}{\mu_{ij^*}}.
\end{align*}
\item All virtual queues $Q_j$ are decreased at a rate $N_j$ per time unit whenever they are positive. The cluster virtual queues $L_j$ for $j \in l$ are decreased when $\sum_{j \in l} \beta L_j > W_l$; if this is the case, we decrease
\[
L_j \mapsto L_j - N_j, ~~ \forall j \in l.
\]
\end{enumerate}
Note that the algorithm is still {\em distributed}, since additional communication needs to happen only between the charging stations within the cluster. Namely, stations must communicate their updated values of $L_j$ whenever these increase, and some entity must decide to decrease the $L_j$ within the cluster once their sum exceeds $\beta^{-1} W_l$.

Similarly to the analysis of GPD, as $\beta \to 0$ the routing pattern produced by the LB algorithm will converge to the optimal solution to \eqref{eqn:LB}. The virtual queues $\beta L_j$ will converge to the optimal dual variables corresponding to the constraints defining the loads, \eqref{load}. Further, under Assumptions~\ref{ass: FCLT} and \ref{ass: CRP}, a variant of Theorem~\ref{thm: FCLT} can be shown to hold, with deviations of the routing process from its nominal value given by some Brownian motion; however, the Brownian motion will now be more cumbersome to define.

In Section~\ref{sec: simulations}, numerical work shows that replacing the GPD algorithm by the LB algorithm can lead to substantial reduction in the queueing delays encountered by the vehicles.

A feature of the GPD and LB algorithms is that they can be slow to reach convergence. Specifically, the typical time scale of GPD and LB is $\beta^{-1}$; so as $\beta \to 0$, the algorithms get more precise but slower at reacting to shocks. A possibility to mitigate this is to run GPD (or LB) in the background to identify the set of edges in the basic activity tree, but then run a different algorithm once that tree has been identified. We will return to this idea in Section~\ref{sec: 0-infty}, where we present one such possible algorithm.
 
\section{The case of No costs}\label{sec: 0-infty}
We now consider the special case of $c_i(j) \in \{0,\infty\}$ for all $i$ and $j$, corresponding to vehicles that are indifferent to the choice of charging station provided they can reach it (and it has compatible technology). This corresponds to the question of optimizing the throughput, i.e. the number of vehicles that are served. In this setting, the parameter $\beta$ in the GPD algorithm is irrelevant, because the decisions will be the same for all values of $\beta$. In particular, we conclude that the GPD algorithm always generates routing rates in the feasible region, provided $\Lambda$ is nonempty.

Note that if the feasible region is nonempty, the optimal dual variables corresponding to the capacity constraints \eqref{capacity} are $q^*_j = 0$. While of course we will not have $Q_j(t) = 0$ at all times for the virtual queues under GPD, they will be close to 0, and will periodically hit 0 unless the system is overloaded.

If the interarrival and service times are exponential, there will be well-defined (unique) steady-state routing rates $\lambda^*_{ij}$, which we know are feasible for \eqref{STABILITY}. However, if $c_i(j) \in \{0, \infty\}$, the optimal solution to \eqref{STABILITY} will certainly \emph{not} be unique. Determining to which optimal solution the routing rates of GPD converge is difficult; this is another reason to use some form of load-balancing.

Note that for the GPD algorithm, in the case of zero costs, convergence time is not an issue: routing rates are ``always'' feasible. However, for the LB one, the quantity $\beta$ makes an appearance, since we decrement the cluster queues $L_j$ when $\sum_{j \in l} \beta L_j = W_l$. Picking values of $\beta$ that are too small will result in imperfect load balancing (although typically will reduce the maximal loads somewhat). A better, and faster, possibility is to use the Freest Charger Shortest Queue (FCSQ) algorithm, modelled on \cite{LQFS}.

\vspace{1cm}\noindent{\bf FCSQ Algorithm:}
\begin{enumerate}
\item Identify the basic activity tree (e.g. by running GPD in the background).
\item When a vehicle of type $i$ requests service, if some charger with $c_i(j) < \infty$ and $\mu_{ij} > 0$ is available, route the vehicle to the charging station with the largest fraction of free chargers.
\item If no charger is available in a charging station, the charging station sends back an estimate of the time when the vehicle will enter queueing. The vehicle joins the charging station where this time is earliest.
\end{enumerate}
The last step of the algorithm is deliberately vague, because many versions should give essentially similar results. In the simulation presented in Section~\ref{sec: simulations}, 
we add up the (future) charging times of all of the EVs in the queue and divide by the number of chargers at the station. We could add to that the residual charging times of the vehicles that already charging; or we could replace the exact times by estimates, reporting $N_j^{-1} (\sum_i \mu_{ij}^{-1} n_i)$, where $n_i$ is the number of queued vehicles of type $i$.

The FCSQ algorithm should guarantee that, when there are free chargers, they are divided fairly among the various charger pools, while when there are queues, they have similar waiting times at all the chargers. This means that when there is a sharp spike in arrivals at one of the charger pools, these arrivals are quickly spread out as far as possible, helping the spike dissipate faster. Strong analytic guarantees are available for a similar algorithm (Longest Queue Freest Server, see \cite{LQFS}) in the case when the non-zero charging rates $\mu_{ij}$ depend only on the charging station technology (i.e. $\mu_{ij} = \mu_j$), and partial results are available for general parameters. Similar algorithms are common in the queueing literature, see for example \cite{Gurvich_Whitt}, or \cite{LAP} and references therein.

\section{Numerical Illustration}\label{sec: simulations}
We now consider the toy network shown in Fig.~\ref{fig: toy example} to illustrate the performance of the various algorithms presented in this work. 
There are two EV classes, three charger classes, and each EV class can use two of the charger classes. The service rates $\mu_{ij}$ are
\[
\mu = \begin{pmatrix}
1 & 3 & 0\\ 0 & 1 & 2
\end{pmatrix}.
\]
We use $N_j = 20$ for all $j$. We simulate 10,000 EV arrivals. The first 5,000 arrivals are generated using arrival rates $\lambda = (2.5,2.2) \times 20$; for the second 5,000 arrivals, we reverse the arrival rates to obtain $\lambda = (2.2, 2.5) \times 20$. We do this to illustrate the effect of a change in the arrival pattern that does not change the basic activity tree
on the algorithms under consideration.
\begin{figure}[htbp]
\centering
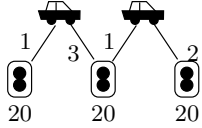
\caption{The diagram for the toy example indicating $N_j$ and $\mu_{ij}$.}
\label{fig: toy example}
\end{figure}
Note that the system is heavily loaded, but not overloaded: the load-balancing solution to the linear program \eqref{eqn:LB}, achieved for sufficiently high values of $W_1$, has $\rho = 0.91$ for the first arrival pattern, and $\rho = 0.97$ for the second arrival pattern. In a conventional single-server queue, such heavy loads would result in long queues and very high probability of having to queue; in the many-charger setting, however, approximately 70\% of the arriving vehicles are taken into service without having to queue at all.

Fig.~\ref{fig: GPD} presents the delays observed when running the GPD algorithm. The three curves present the delays encountered by the vehicles routed to each of the three charging stations. (One of the curves stays at 0 throughout the simulation.) The vertical line indicates the time when the arrival pattern changes. We see that before the change in the arrival pattern, a large queue had been forming at station 2 (dashed line); after the change, it slowly goes away, but a sizeable queue forms at station 3 (dotted line) instead.
\begin{figure}[htbp]
\centering
\includegraphics[width=10cm]{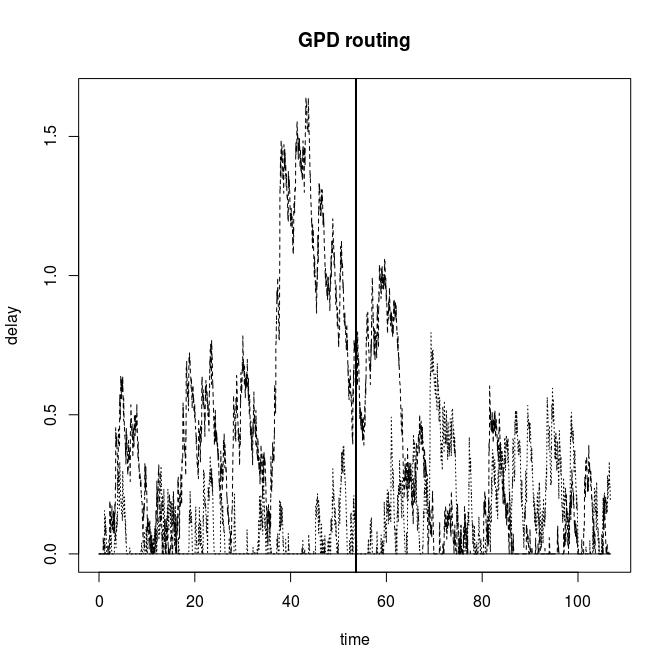}
\caption{Delays produced by GPD routing}
\label{fig: GPD}
\end{figure}

In Fig.~\ref{fig: LB}, we simulate the LB algorithm, putting all of the chargers into a single cluster. We use $W_1 / \beta = 1000$. The scale is the same as for the GPD routing; so we see that the delays have gotten shorter. (The maximal delay has dropped from 1.64 to 1.3.) We also see the marked change in the pattern after the change in the arrival rate pattern: the delay at station 2 (dashed line) drops nearly to zero, however, it takes a while for that to occur. The behavior of the LB algorithm would be improved if we considered a larger system ($N_j > 20$).
\begin{figure}[htbp]
\centering
\includegraphics[width=10cm]{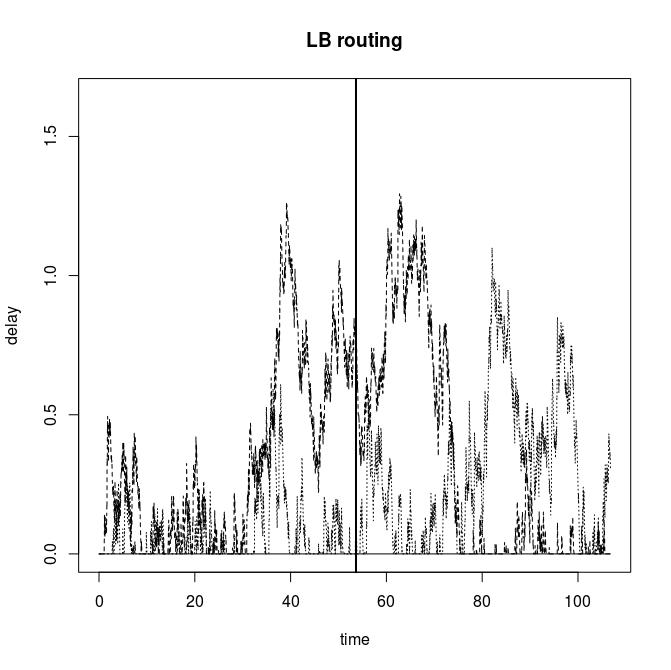}
\caption{Delays produced by LB routing}
\label{fig: LB}
\end{figure}

Fig.~\ref{fig: FCSQ} shows the delays arisint from FCSQ routing on the same arrival data. Note that the largest delay is now only 0.94, more than a third smaller than under GPD. Also, when there is queueing delay, it is nearly the same at all stations, meaning that there is little incentive for any one vehicle to disobey the algorithm and go elsewhere. This is in shaprp contrast to GPD, which tends to produce highly unequal queueing delays.
\begin{figure}[htbp]
\centering
\includegraphics[width=10cm]{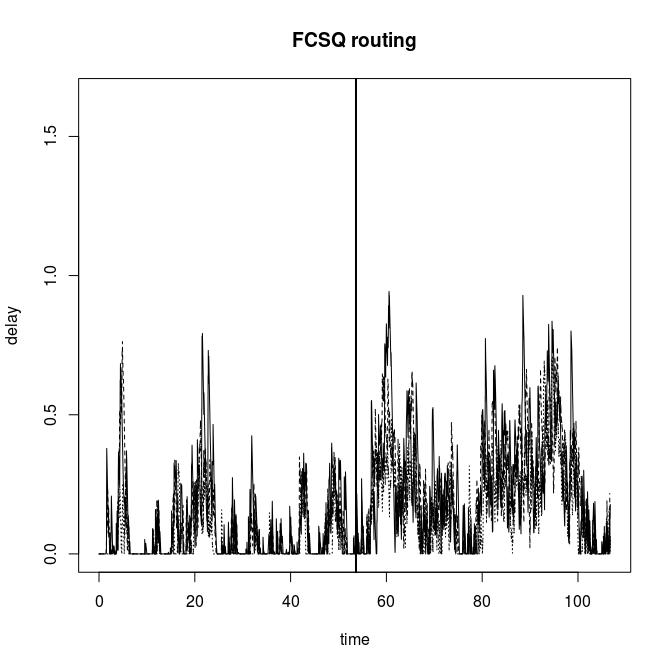}
\caption{Delays produced by FCSQ routing}
\label{fig: FCSQ}
\end{figure}

\section{Discussion and future work}\label{sec: discussion}
We have constructed a set of algorithms -- GPD, LB, and FCSQ -- which collectively route the EVs to chargers in the network in such a way as to avoid long queueing and adjust to changing demand patterns whenever this is possible. The LB algorithm s designed to balance the loads across different charging stations; the FCSQ algorithm 
corresponds to a faster version of LB, which, however, can only be run after LB or GPD has been used to identify the basic activities along which routing should happen. All of the algorithms are distributed, and therefore can scale to a large network.

The algorithms have implications for building the charging network. As we observed in Section~\ref{Halfin-Whitt}, square-root overstaffing is needed for the operation of a 
system with many chargers in each pool and a small value of the parameter $\beta$ with short queueing delays. Moreover, the overstaffing is aggregated over all of the charging stations in a single connected component of the basic activity graph. Therefore, for greater freedom in the placement of extra charging facilities it is desirable to keep these connected components large. That is, in contrast to the facility locations problem, we not only allow, but {\em encourage} EVs of a single type to visit multiple charger types.

Two important extensions of this work are of interest. First, it is desirable to arrange charging stations so that the basic activity graph for the system will have large connected components. However, this is a highly nonlinear and nonconvex requirement. We would like to create an algorithm for generating arrangements of charging stations satisfying this condition.
Second, we would like to incorporate the effect of batteries onto this system. Because of the high load on the electrical network produced by charging an electric vehicle, it is desirable to add battery capacity to the charging stations to smooth the peak demand. This introduces an additional algorithmic challenge by introducing demand mobility in time in addition to in space.

Finally, we remark that due to technical challenges, some of the results concerning the Halfin-Whitt regime remain conjectural. Specifically, there is a gap between the finite-time-horizon results on the Brownian motion approximations, and steady-state quantities such as the probability of an arriving request being asked to queue. However, 
this gap is not very important in practice because of the diurnal demand pattern.

\appendix[Proof of Theorem~\ref{thm: FCLT}]\label{sec: proof of thm}
The overall technique of proof of Theorem~\ref{thm: FCLT} follows the proof of \cite[Theorem 6.3]{Shadow_routing}. We outline both the argument and the necessary changes below. 

\begin{proof}[Sketch of proof of Theorem~\ref{thm: FCLT}]
To show the first assertion of the theorem, we begin by demonstrating that if the virtual queues are initially close to their  equilibrium values, then the deviations must remain small throughout the entire time interval $[0,t]$. The version of this argument in \cite[Theorem 6.1]{Shadow_routing} does not hold in our case; however, the following result is true.

\begin{lemma}[Version of Theorem 6.1 of \cite{Shadow_routing}]\label{lm: thm 6.1}
If $Q^r_j(0) - (\beta^r)^{-1} q^*_j = o(\sqrt{r})$, then uniformly on any compact set $t \in [0,T]$ we have
\[
\beta^r Q^r_j(t) - q^*_j = o(1), ~~ \forall j,  ~~\text{uniformly on any compact set $t \in [0,T]$.}
\]
Consequently,
\[
\left(c_i(j) + \frac{\beta^r Q^r_j(t)}{\mu_{ij}}\right) - \left(c_i(j') + \frac{\beta^r Q^r_{j'}(t)}{\mu_{ij'}}\right)  = o(\sqrt{r}), ~~ \forall (ij), (ij') \in \CE,
\]
uniformly on any compact set $t \in [0,T]$.
\end{lemma}
We postpone the proof of the lemma to proceed with the rest of the argument. The first claim of the lemma implies that, for $r$ sufficiently large, vehicles will not be routed along non-basic activities, because differences in virtual queue sizes dominate the fluctuations in the size of each virtual queue. The second claim then determines the proportion of vehicles routed along each of the basic activities. Because the differences between the quantities controlling the routing are much smaller than the size of the Brownian fluctuations in the arrivals, the Brownian fluctuations are unaffected by them.

The stationary result is proved by considering the fluid limits of the form $\sqrt{r}(Q^r(\sqrt{r}u) - q^*_j)$. The analysis there shows that in steady state, this deviation is $o(\sqrt{r})$; indeed, it should be $o((\beta^r)^{-1})$. We refer the reader to \cite[Theorem 6.3, 6.4]{Shadow_routing} for details.
\end{proof}

We now outline the proof of Lemma~\ref{lm: thm 6.1}.
\begin{proof}[Sketch of proof of Lemma~\ref{lm: thm 6.1}]
Both claims of the lemma are proved using the technique of local fluid limits. For details, the reader may consult e.g. \cite[Section 8]{Mandelbaum_Stolyar}.

For the first assertion, define for all $j$ the local-fluid-scaled queueing process
\[
q^{(r,m)}_j(u) = \beta^r Q^r_j(t^m + \beta^r u) - q^*_j, ~~ u \in [0,T], ~~ m = 1, \dotsc, r \beta^r.
\]
It is standard to show that, as $r \to \infty$, each sequence $q^{(r,m)}_j(\cdot)$ converges to the family of Lipschitz processes $q^m_j(\cdot)$ satisfying certain differential equations whenever they are defined. Our goal is to show that $\max_j q^{(r,m)}_j(u)$ decreases whenever it is positive, and does so at some rate bounded away from 0. Note first that for finite values of $q^{(r,m)}$, the values of $Q^r$ are still almost proportional to their nominal values $q^*_j$, and therefore all requests must be routed along basic activities only. Moreover, consider the set of charger types that have the largest value of $q^m(u)$. The routing rule of GPD algorithm ensures that these charger types will have only those vehicles routed to them that cannot be routed elsewhere. By Assumption~\ref{ass: CRP}, if the maximal value of $q^m$ is positive, then the virtual queues for which $q^m$ is maximal will have arrivals that are smaller than nominal, and of course will be decreased at the nominal rate. Consequently, when the maximal value of $q^m$ is positive, it must decrease towards its equilibrium value of 0. The time $T$ is picked to be sufficiently large so that, e.g., if $q^m_j(0) \leq 1$ for all $j$, then $q^m_j(T) = 0$ for all $j$. We obtain that on a single time interval, $q^{(r,m)}(u)$ is very likely to stay close to 0. Finally, we show that this holds on all $r \beta^r$ subintervals simultaneously. Suppose not, and pick the subsequence of intervals $m(r)$ along which $q^{(r,m)}$ first crosses level $\epsilon < 1$. Pick a subsequence along which $q^{(r,m)}(0)$ converges, and observe that the corresponding local fluid limit contradicts our previous assertions.

For the second assertion, we define the slightly different local-fluid-scaled processes
\[
\hat{q}^{(r,m)}_j(u) = r^{-1/2} \left(\beta^r Q^r_j(t^m + r^{-1/2} u) - q^*_j\right), ~~ u \in [0,T], ~~ m = 1, \dotsc, \sqrt{r}.
\]
Very similar arguments tell us that for finite values of $\hat{q}^{(r,m)}$, the quantities $c_i(j) + \mu_{ij}^{-1} \beta_j Q^r_j$ are approximately proportional to their nominal values, and therefore routing decisions are very close to nominal. This easily implies that the quantities $c_i(j) + \mu_{ij}^{-1} \beta_j Q^r_j$ stay close to their nominal values.
\end{proof}

\bibliographystyle{IEEEtran}
\bibliography{IEEEabrv,EV}

\end{document}